\documentclass{eurocg13}

\usepackage{tikz}
\usetikzlibrary{positioning}
\usetikzlibrary{calc}
\usetikzlibrary{intersections}
\usetikzlibrary{arrows}

\definecolor{darkgreen}{RGB}{0,153,0}
\definecolor{darkpurple}{RGB}{178,122,255}
\definecolor{darkorange}{RGB}{255,128,0}

\usepackage{wrapfig}

\usepackage[hidelinks]{hyperref}
\usepackage{url}

\title{\vspace{-2ex}Constructing Complicated Spheres\\[.5ex]\small (Extended Abstract)}

\author{Mimi Tsuruga \hspace{4em} Frank H. Lutz\\
Technische Universit\"at Berlin\\
\tt	\{tsuruga,lutz\}@math.tu-berlin.de\vspace{-3ex}}

\begin{document}
\maketitle

\abstract{
Fast and efficient homology algorithms are in demand in the applied sciences for analyzing solid materials and proteins, processing digital imaging data, or pattern classification among others. Recent advances employ discrete Morse theory as a preprocessor. Research in this area has lead to the need to find complicated test examples. We present an infinite series of examples that have been constructed to test some of the latest algorithms under development. This family of 4-spheres (known as the Akbulut--Kirby spheres) is based on a handlebody construction via finitely presented groups.
}

\section{Motivation}
\subsection{Homology Algorithms}

Computing the homology of simplicial or more general cell complexes can quickly become too cumbersome by hand. Today, however, computers can be used to analyze even amazingly large complexes~\cite{redhom, chomp}. Speeding up homology calculations is now an active field of research in computational topology.
    
Computing homology involves the computation of the Smith normal form of matrices that encode geometric information about the complex. Current algorithms for computing the Smith normal form, though polynomial~\cite{kannan}, are too slow for explicit computations on large examples. However, computation time can be dramatically improved by employing a topological preprocessor.

Discrete Morse theory provides an excellent foundation for reducing a large complex to a simplified complex of the same topological (PL) type~\cite{Joswig2004pre,KaczynskiMischaikowMrozek2004}. We can then replace the matrix for which we want to find the Smith normal form with a much smaller one. The problem of finding an optimal discrete Morse function, however, is $\mathcal{NP}$-hard~\cite{joswigpfetsch,lewiner}.

Benedetti and Lutz~\cite{benedettilutz} propose using random discrete Morse theory to search for discrete Morse functions with few critical cells. Surprisingly, this approach frequently finds the optimum even for large examples. In fact, it has been difficult to provide complicated examples on which the random heuristics performs poorly.

The aim of this article is to provide an infinite series of such examples.
    
\subsection{Complicated Spheres}

Simplicial complexes can pose challenges for discrete Morse algorithms for various reasons. For example, if the underlying topological space is complicated or if the triangulation itself is complicated. Here we look at the latter case and construct a series of triangulations of the 4-dimensional sphere $S^4$, which has trivial topology.\footnote{We study spheres because spheres can admit perfect Morse functions, i.e., with only two critical cells.} The examples in the series, however, are composed in a rather non-trivial manner.

\section{Akbulut--Kirby Spheres}

\subsection{Background}\label{background}
Topologists in recent years have been working on finding counterexamples for the smooth Poincar\'e conjecture in dimension 4. These so-called exotic spheres are smooth manifolds that are homeomorphic, but not diffeomorphic to the standard sphere $S^4$, i.e., the round sphere with the usual differentiable structure. Various candidates have been introduced, many of which were later shown to be standard.

We have built explicit triangulations of one such candidate series proposed by Akbulut and Kirby~\cite{potentialcounterexample}, which goes back to an example by Cappell and Shaneson. Gompf~\cite{killingAK} showed that one sphere from the series is standard, i.e., diffeomorphic to $S^4$. Akbulut~\cite{akbulut} later proved that all examples in the whole infinite series are standard.

In the following section, we will give a detailed description of our construction of triangulations of the Akbulut--Kirby spheres.

\subsection{Construction}

Before going into the details of our construction of triangulations for the Akbulut--Kirby spheres, let us start with a quick overview of what's to come. The main idea for the spheres begins with finitely presented groups 
\begin{align*}
 G = \left< g_1, \dots, g_l \mid r_1, \dots, r_m \right>,
\end{align*}
where $g_i$ are generators and $r_i$ are relators in the $g_i$'s. It is undecidable whether a given finitely presented group is the trivial group~\cite{novikov}. In fact, there are even group presentations with only two generators and two relators that are not obviously trivial.

Let $G = \left< x, y \mid xyx = yxy ; x^r = y^{r-1}\right>$. To see that this group is trivial, rearrange the first relator to get
\begin{align*}
 y &= x^{-1} y^{-1} x y x.
 \end{align*} 
Then $y^r =x^{-1} y^{-1} x^r y x = x^{-1} y^{r-1} x = x^r$ by substituting the second relator twice. This yields $y^r=y^{r-1}$ or $y=e$ and hence $x=e$. Thus $G$ is the trivial group.

Next, take a 5-ball and attach two 1-handles representing each of the generators $x$ and $y$. Then glue in two 2-handles along thickened curves representing the relators $xyx=yxy, x^r=y^{r-1}$. The resulting object is a topological 5-ball. 

Gompf~\cite{killingAK} proved that for all the presentations $G=\left<x,y \mid xyx=yxy,x^r=y^{r-1}\right>$ the resulting boundary 4-spheres are, in fact, standard 4-spheres.

We built simplicial complexes following this recipe. But we begin our construction in dimension~3.

Start with a 3-ball and attach two 1-handles. Then choose two curves in the resulting handlebody that represent the relators. Next we take the product of the 3-dimensional handlebody $H^3$ with an interval $I$ twice. We then find (copies of) the two curves in the boundary of the 5-dimensional handlebody $H^3 \times I \times I$.

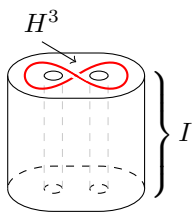
\begin{wrapfigure}{l}{.18\textwidth}
\vspace{-8ex}%
\begin{center}

\begin{tikzpicture}[scale=.3]

\draw ($(90:2 and 1)+(0,5)$) arc (90:270:2 and 1)
      -- ++(2,0)
      arc (270:450:2 and 1)
      -- cycle;

\draw (180:2 and 1) arc (180:270:2 and 1)
      -- ++(2,0)
      arc (270:360:2 and 1);
\draw[dashed] ($(0:2 and 1)+(2,0)$) arc (0:90:2 and 1)
      -- ++ (-2,0)
      arc (90:180:2 and 1);

\draw (4,0) -- (4,5) (-2,0) -- (-2,5);

\draw[dashed] (0,0) ellipse (.4 and .2);
\draw[dashed] (2,0) ellipse (.4 and .2);
\draw[dashed,black!20] (-0.4,0)--(-0.4,5)  (0.4,0)--(0.4,5)
              (1.6,0)--(1.6,5)  (2.4,0)--(2.4,5);

\draw (0,5) ellipse (.4 and .2);
\draw (2,5) ellipse (.4 and .2);

\node (par) at (5.1,2.5) {$\left.\rule{0cm}{1cm} \right\} I$}; 

\draw[->] (-.5,6.5) node[above] {$H^3$} -- (1,5.6);


%
 

\draw[red,thick] (1,5) .. controls (-2,7) and (-2,3) ..
                 (1,5) .. controls (4,7) and (4,3) ..
                 (1,5);

\draw[white,line width=1pt](1.5,4.82) -- (.5,5.44);
\draw[white,line width=1pt](1.5,4.58) -- (.5,5.16);

\end{tikzpicture}
 \end{center}
\vspace{-6ex}
 \caption{\small{The red curve in the interior of $H^3$ is on the boundary of $H^3 \times I$.}}
 \label{crossi}
 \vspace{-2ex}
\end{wrapfigure}

Consider Fig.~\ref{crossi} in which the red curve lies in the interior of the handlebody $H^3$. After taking 
 the product $H^3 \times I$, we find (a copy of) the red curve on the boundary of $H^3 \times I$. On the boundary of $H^3 \times I \times I$ we thicken each of the two curves representing the generators $xyx=yxy$, $x^r = y^{r-1}$ to a solid 4-torus along which the two 2-handles are glued in. While the two curves are twisted in $H^3$, the two solid 4-tori untwist in the boundary of $H^3 \times I \times I$. 

By gluing in the two 2-handles, we end up with a topological 5-ball whose boundary (if smoothed appropriately) is a homotopy 4-sphere~\cite{potentialcounterexample}. For 4-manifolds, the categories DIFF (differentiable) and PL (piecewise linear) coincide, so our triangulations provide model spaces for the smooth setting.

\subsubsection{The Space.}

\begin{figure}[ht]
\begin{center}

\begin{tikzpicture}[scale=.3]

\fill[darkorange!50] (180:10 and 3) -- ($(180:10 and 3)+(0,-1)$)
       arc (180:360:10 and 3) --(360:10 and 3);

\draw ($(180:10 and 3)+(0,-1)$) arc (180:360:10 and 3);
\draw (180:10 and 3) -- ($(180:10 and 3)+(0,-1)$)
      (360:10 and 3) -- ($(360:10 and 3)+(0,-1)$);
      
\filldraw[fill=darkorange,draw=black] (0,0) ellipse (10 and 3);
\filldraw[fill=darkpurple,draw=black] (0,5) ellipse (10 and 3);

\filldraw[fill=darkgreen,draw=black] (-10,0) rectangle (10,5);

\filldraw[fill=white,draw=black] (0,0) ellipse (3 and 1);
\filldraw[fill=white,draw=black] (0,5) ellipse (3 and 1);

\begin{scope}
 \clip (0,0) ellipse (3 and 1);
 \filldraw[fill=darkgreen!50,draw=black]
 ($(180:3 and 1)+(0,-1)$) arc (180:0:3 and 1) -- (0:3 and 1)
 arc (0:180:3 and 1);
\end{scope}

\begin{scope}
 \clip (0,5) ellipse (3 and 1);
 \filldraw[fill=darkpurple!50,draw=black]
 ($(180:3 and 1)+(0,4)$) arc (180:0:3 and 1) -- ($(0:3 and 1)+(0,5)$)
 arc (0:180:3 and 1);
\end{scope}

 \node (x) at (0,7) {$x$};
 \node (y) at (0,-2) {$y$};
\end{tikzpicture}
\end{center}
\vspace{-4ex}
 \caption{The handlebody $H^3$ consisting of a 3-ball (in green) with two 1-handles (purple for the generator $x$ and orange for the generator $y$).}
 \vspace{-2.5ex}
 \label{body}
\end{figure}
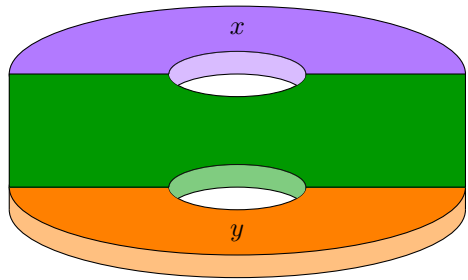

We begin in dimension 3. In particular, we want to have a 3-ball with two 1-handles which represent our two generators $x$ and $y$, see Fig.~\ref{body}.

Notice that we have used three colors to distinguish the different parts of our space. Keep this picture in mind when we move on to the next steps. The green section in the center is our 3-ball. To this 3-ball we attach one 1-handle which we label $x$ (in purple) and another 1-handle which we label $y$ (in orange). 

We will choose two curves that represent the two relators inside this space.

\subsubsection{Planning the Curves.}

\begin{figure}[h]
\begin{center}
\begin{tikzpicture}
    \node[anchor=south west,inner sep=0] at (0,0) {\includegraphics[width=.4\textwidth]{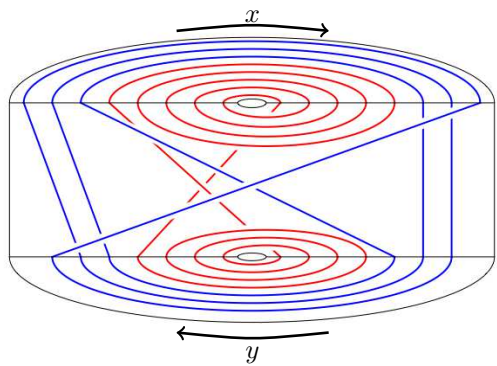}};
       
    
    \coordinate (top) at (3.45,4.15);
    \coordinate (btm) at (3.45,0);
    \draw[->,line width=1pt]  ($(top)+(-1,-.05)$) .. controls ($(top)+(0,.05)$) .. ($(top)+(1,-.05)$);
    \node (x) at ($(top) + (0,0.15)$) {$x$};
    
    \draw[->,line width=1pt]  ($(btm) + (1,.1)$) .. controls (btm) .. ($(btm) + (-1,.1)$);
    \node (y) at ($(btm) + (0,-.2)$) {$y$};
\end{tikzpicture}
\end{center}
\vspace{-4ex}
\caption{The two curves of the Akbulut--Kirby sphere $xyx=yxy$ (blue) and $x^5=y^4$ (red) in dimension 3.} \label{spaghetti}
\end{figure}

We draw the two curves which represent the relators. The two curves are $xyx=yxy$, which we call the blue curve, and $x^r=y^{r-1}$, which we call the red curve (see Fig.~\ref{spaghetti} for the Akbulut--Kirby sphere with $r=5$).

\vspace{-2ex}%

The crossings indicate which sections of the curves run over or under. The embedding of the curves can be arbitrarily chosen as long as they run over the two 1-handles $x$ and $y$ as specified by the relators. When we move up in dimension in Step~6, the thickened curves will untangle no matter how twisted or knotted they were in dimension~3. The particular arrangement shown in Fig.~\ref{spaghetti} was chosen to simplify the construction in Step~5.

To understand Fig.~\ref{spaghetti}, begin with the blue curve read as $xyxy^{-1}x^{-1}y^{-1}=e$. Remember that we are on a 3-ball with two 1-handles $x$ and $y$. These handles are oriented as indicated by the black arrows.

For the blue curve, we begin at the blue crossing in the center of Fig.~\ref{spaghetti}. Follow the bottom curve at the crossing and go up (and left) towards the $x$-handle. Then along $x$, back through the ball, then $y$, then $x$, then all the way across to the far bottom left, along $y$ in the reverse direction, then $x$ in reverse, then $y$ in reverse and then finally return to our starting point to close the curve. The red curve $x^5y^{-4}=e$ is represented similarly.

\subsubsection{Thickening of the Curves to Solid Tori.}

Fig.~\ref{spaghetti} displays the blue and red curves in $H^3$ which will be thickened to solid 3-tori. We choose triangular prisms to compose the solid 3-tori.%
\vspace{-2ex}%
\begin{figure}[h]
 \begin{center}

\begin{tikzpicture}[scale=.4]

\begin{scope}
 \clip (1.5,0) -- (1,2) -- (7,2) -- (7.5,0) -- cycle;
 \fill[blue,opacity=.5] (1,0) rectangle (7.5,2);
\end{scope}

\fill[blue,opacity=.5] (1.5,0) -- (1,2) -- (0,.5) -- cycle;

\begin{scope}
 \clip (13.5,0) -- (13,2) -- (16,2) -- (16.5,0) -- cycle;
 \fill[blue,opacity=.5] (13,0) rectangle (16.5,2);
\end{scope}

\fill[blue,opacity=.5] (13.5,0) -- (13,2) -- (12,.5) -- cycle;

 \draw[blue,thick] (1.5,0) -- (1,2) (1,2) -- (0,.5) (0,.5) -- (1.5,0);
 \draw[blue,thick] (4.5,0) -- (4,2);
 \draw[dashed,blue!50!gray] (4,2) -- (3,.5) -- (4.5,0); 

\begin{scope}[shift={+(3,0)}]
 \draw[blue,thick] (4.5,0) -- (4,2);
 \draw[dashed,blue!50!gray] (4,2) -- (3,.5) -- (4.5,0); 
\end{scope} 
 
\draw[blue,thick](1,2) -- (7,2)    (13,2) -- (16,2) 
		(1.5,0) -- (7.5,0) (13.5,0) -- (16.5,0);
\draw[dashed,blue!50!gray] (0,.5) -- (6,.5)   (12,.5) -- (15,.5);
 
 \fill (8.25,1) circle (2pt);
 \fill (9,1) circle (2pt);
 \fill (9.75,1) circle (2pt); 
 
\begin{scope}[shift={+(12,0)}]
 \draw[blue,thick] (1.5,0) -- (1,2) (1,2) -- (0,.5) (0,.5) -- (1.5,0);
 \draw[blue,thick] (4.5,0) -- (4,2);
 \draw[dashed,blue!50!gray] (4,2) -- (3,.5) -- (4.5,0); 
\end{scope}

 \node[above left] (v1) at (0,.5) {$1$};
 \node[below] (v2) at (1.5,0) {$2$};
 \node[above] (v3) at (1,2) {$3$};
 \node[above left] (v4) at (3,.5) {$4$};
 \node[below] (v5) at (4.5,0) {$5$};
 \node[above] (v6) at (4,2) {$6$};
 \node[above left] (v7) at (6,.5) {$7$};
 \node[below] (v8) at (7.5,0) {$8$};
 \node[above] (v9) at (7,2) {$9$};
 \node[above left] (v34) at (12,.5) {$34$};
 \node[below] (v35) at (13.5,0) {$35$};
 \node[above] (v36) at (13,2) {$36$};
 \node[above left] (v37) at (15,.5) {$1$};
 \node[below] (v38) at (16.5,0) {$2$};
 \node[above] (v39) at (16,2) {$3$};

\end{tikzpicture}
 \end{center}
 \vspace{-4ex}
\caption{The blue curve in Fig.~\ref{spaghetti} is thickened to a chain of triangular prisms.}
\vspace{-1ex}
\label{triprisms}
\end{figure}
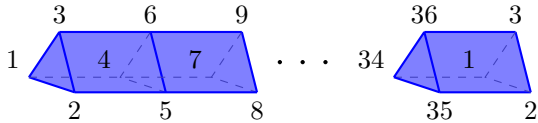
In Fig.~\ref{triprisms}, we have the prisms for the blue curve for the Akbulut--Kirby sphere with $r=5$ (Fig.~\ref{spaghetti}). The vertex labels repeat at the ends as they are identified to form a loop.

Notice that the blue curve in Fig.~\ref{spaghetti} can be cut up into 12 strands (3 on the $x$-handle, 3 on the $y$-handle, and 6 inside the middle ball) every time the curve crosses into a different section of the handlebody (recall Fig.~\ref{body}). One triangular prism is used for each of these strands. So we used ($12 \times 3 =$)36 vertices to build the blue solid torus as indicated in Fig.~\ref{triprisms}.

We can systematically and coherently break down these prisms into tetrahedra by using product triangulations as described in~\cite{lutztriangulate} and references therein. This method guarantees that neighboring simplices match up nicely. For example, we do not want to introduce different diagonals for a rectangular face of one of the prisms.

\subsubsection{Fill the Handles.}

The construction of the 1-handles is simple since they only have triangular prisms that run parallel along each other. However, we want to avoid having any unwanted identifications along the handles. So we add buffers between neighboring triangular prisms using rectangular prisms.%
\vspace{-2ex}
\begin{figure}[ht]
\begin{center}
\begin{tikzpicture}[scale=.3]


\begin{scope}
\clip (16,0) -- ($(15,0)+(60:1)$) arc (12:170:2.3 and 1.5) -- (11,0); 
\fill[red] (11,0) rectangle (16,2.5);
\end{scope}
\fill[white] (15,0) arc (12:170:1.8 and 1.2);

\draw ($(25.5,0)+(60:1)$) arc (12:170:13.02 and 6.5);
\draw ($(24,0)+(60:1)$) arc (12:170:11.5 and 6);
\draw ($(22.5,0)+(60:1)$) arc (12:170:9.98 and 5.5);
\draw ($(21,0)+(60:1)$) arc (12:170:8.43 and 5);
\draw ($(19.5,0)+(60:1)$) arc (12:170:6.9 and 4.5);
\draw ($(18,0)+(60:1)$) arc (12:170:5.37 and 3.5);
\draw ($(16.5,0)+(60:1)$) arc (12:170:3.83 and 2.5);
\draw ($(15,0)+(60:1)$) arc (12:170:2.3 and 1.5);
\draw (15,0) arc (12:170:1.8 and 1.2);

\draw ($(0,0)+(60:1)$) -- ($(10.5,0)+(60:1)$)
           ($(15,0)+(60:1)$) -- ($(25.5,0)+(60:1)$);
\draw (0,0) -- (11,0) (15,0) -- (26,0);

\foreach \i in {1,2, ...,3}
\filldraw[fill=blue,draw=black] ($(\i,0)-(1,0) + {\i - 1}*(.5,0)$) -- ++(0:1) -- ++(120:1) -- cycle;

\foreach \i in {4,5, ...,8}
\filldraw[fill=red,draw=black] ($(\i,0)-(1,0) + {\i - 1}*(.5,0)$) -- ++(0:1) -- ++(120:1) -- cycle;

\foreach \i in {11,12, ...,15}
\filldraw[fill=red,draw=black] ($(\i,0)-(1,0) + {\i - 1}*(.5,0)$) -- ++(0:1) -- ++(120:1) -- cycle;

\foreach \i in {16,17, ...,18}
\filldraw[fill=blue,draw=black] ($(\i,0)-(1,0) + {\i - 1}*(.5,0)$) -- ++(0:1) -- ++(120:1) -- cycle;

\end{tikzpicture}
\end{center}
\vspace{-4ex}
\caption{Filling the $x$-handle of Fig.~\ref{spaghetti}.\vspace{-2ex}
} \label{handlefill}
\end{figure}
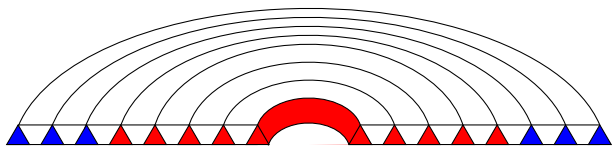
Fig.~\ref{handlefill} gives a front view of the $x$-handle. The three blue triangles on the two ends are part of the blue solid 3-torus running over the $x$-handle in Fig.~\ref{spaghetti}. The rectangular prisms between the triangular prisms are white. We need 7 rectangular buffer prisms for the $x$-handle and 6 rectangular buffer prisms for the $y$-handle.

\subsubsection{Fill the Ball.}

\begin{figure}[h]
\begin{center}
 \includegraphics[width=0.47\textwidth]{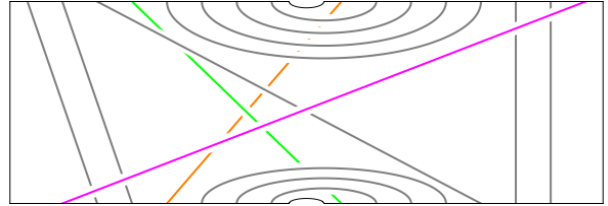}
\end{center}
\vspace{-15pt}
\caption{The portion that will be filled to become the 3-ball. The three colored lines go above or below.
} \label{floors}
\end{figure}

First we explain the reason for choosing the particular arrangement of the curves in Fig.~\ref{spaghetti}. Let's look only at the rectangular area representing what will become our 3-ball, see Fig.~\ref{floors}.

\vspace{-1ex}

The (curves that represent) prisms in the middle ball section are not all parallel and some cross each other. The gray lines have no crossings and we place them in, say, Level 0. We then let the pink line go above (Level +1), green goes below (Level -1) and orange further below (Level -2).

To fill the floors, we first glue in rectangular buffer prisms between the parallel strands of the ground floor. Then we glue in 2-dimensional (triangulated) membranes to close the holes between the remaining strands of the ground floor, then two additional membranes to close an upper cupula including the pink strand and further membranes to obtain another two cavities for the two lower floors. The three cavities are then closed by filling in three cones.

The main part of our construction is now complete.

\subsubsection{Go Up In Dimension.}

Having composed $H^3$, we can now take a direct product of it with the unit interval $I=[0,1]$ once to go to 4-dimensions. In this step, we thicken the 3-dimensional solid tori to 4-dimensional solid tori. We then take the direct product again to go to 5 dimensions. This time, we keep the solid 4-tori in one of the two boundary components. The product triangulation procedure~\cite{lutztriangulate} is used each time.

\subsubsection{Glue In the 2-Handles.}

We glue in the two 2-handles along the two solid 4-tori in a canonical way by using four additional vertices each.

Now we have a 5-dimensional ball.

\subsubsection{Take Its Boundary.}

Finally, we take the boundary to obtain our (homotopy) 4-sphere. Notice that all the vertices have been pushed out onto the boundary at Step 6.

The final steps use techniques that are mostly canonical and will be described in more detail in a full version of this article (in preparation).

\section{Results}

\begin{theorem}
The series of Akbulut--Kirby spheres can be triangulated with face vectors
$f=(176$~$+$~$64 \, r$, $2390$~$+$~$1120\,r$, $7820$~$+$~$3840\,r$, $9340$~$+$~$4640\,r$, $3736$~$+$~$1856\,r)$ for $r\ge3$, where the triangulations respect the defining handlebody decompositions as described in the above construction.
\end{theorem}

As we hoped, these triangulated spheres have shown to be substantially complicated. For example for $r=5$, the best discrete Morse vector found, $( 1, 2, 4, 2, 1)$, was reached only 5 times out of 1,000 random runs. On average, we certainly need more than just two critical cells; see Table~\ref{critcells}.
\vspace{-10pt}
\begin{table}[h]
\begin{center}
\begin{tabular}{p{2em} r@{.}l p{1em} p{2em} r@{.}l}
 $r$ & \multicolumn{2}{c}{$c_{\tau}$} &{ }& $r$ & \multicolumn{2}{c}{$c_{\tau}$} \\ \cline{1-3} \cline{5-7}
 \rule{0pt}{12pt}5&29&26 	&{ }& 8&40&07\\
 6 &33&124 	&{ }& 9&42&432\\
 7&37&026	&{ }& 10&46&946
\end{tabular}
\end{center}
\vspace{-18pt}
\caption{Average number of critical cells $c_{\tau}$ needed for $r=5, \dots, 10$ in 1,000 runs.}
\vspace{-2ex}
\label{critcells}
\end{table}

It is important to note that there does not exist an algorithm that can determine whether a (triangulated) manifold of dimension $d>4$ is a sphere and existence of such an algorithm is still open for $d=4$. One heuristic way, however, is to perform (random) bistellar flips, or Pachner moves~\cite{bjornerlutz,pachner} to show that a given triangulation is PL-equivalent to the boundary of a $(d+1)$-simplex. We have shown using \texttt{polymake}~\cite{polymake} that at least for $r=3$ the corresponding Akbulut--Kirby sphere is standard thus reconfirming the statement in this case experimentally.

Another remarkable result of our experiments is that the initial nontrivial group presentation can show up when we compute the fundamental group as the edge-path group according to Seifert and Threlfall~\cite{seifert} and then simplifying the presentation using \texttt{GAP}~\cite{gap}. In fact, we found that the fundamental group of the complex for $r=5$ after some bistellar and (FP) group simplification had two generators $x$ and $y$ and two relators $xyx=yxy$ and $x^5=y^4$. This is exactly the presentation with which we started!

These spheres have recently been used as test examples for the Perseus homology algorithm by Mischaikow and Nanda~\cite{perseus}. The techniques designed to generate these spheres will also be used to construct other PL versions of topologically interesting objects.

\section*{Acknowledgments}
Many thanks to John M. Sullivan for helpful discussions and remarks. We are grateful to Joel Hass and to William~P.\ Thurston for directing us to the Akublut--Kirby spheres as non-trivial constructions of 4-spheres.

{\small \bibliography{.}}

\end{document}